\documentclass[onecollarge,final]{svjour2}
\smartqed  % flush right qed marks, e.g. at end of proof
\usepackage{graphicx}
\usepackage{amsmath}
\usepackage{amsfonts}
\usepackage{float}

\journalname{Archive for Rational Mechanics and Analysis}
\begin{document}

\title{Creating materials in which heat propagates along a line: theory and numerical results
}


\author{Alexander G. Ramm \and Cong Tuan Son Van
}


\institute{Alexander G. Ramm \at
              Mathematics Department, Kansas State University, Manhattan, KS 66506, USA  \\
              \email{ramm@math.ksu.edu}           %  \\
           \and
           Cong Tuan Son Van \at
              Mathematics Department, Kansas State University, Manhattan, KS 66506, USA  \\
              \email{congvan@math.ksu.edu}           %  \\
}

\date{\notused}
% The correct date

\maketitle

\begin{abstract}
In this paper the theory is developed for creating a material in which the heat is transmitted along a given line. This gives a possibility to transfer information using heat signals. This seems to be a novel idea. The technical part of the theory is the construction of the potential $q(x)$. This potential describes the heat equation $u_t = \Delta u - q(x)u$ in the limiting medium which is obtained after the small impedance particles are distributed in a given domain. A numerical method is also established to construct numerically such a potential.
\end{abstract}

\keywords{materials; heat transfer; inverse problems; numerical mathematics.}

\section{Introduction}\label{Intro}
To create materials in which heat propagates along a line, one needs to create a medium in which the heat transfer is governed by the equation
\begin{equation}\label{eqi1}
u_t = \Delta u - q(x)u \quad \text{in } D, \quad u|_S = 0, \quad u|_{t=0}=f(x),
\end{equation}
where $D$ is a bounded domain with a piecewise-smooth boundary $S$, $D = D_0 \times [0, L]$, $D_0 \subset \mathbb{R}^2$ is  a smooth domain orthogonal to the axis $x_1$, $x = (x_1, x_2, x_3)$, $x_2, x_3 \in D_0$, $0 \leq x_1 \leq L$.\\
Such a medium is created by embedding many small impedance particles $D_m$, $1 \leq m \leq M$, into a given domain $D$ filled with a homogeneous material. Let us assume that the distribution of small particles is:
\begin{equation}\label{eqi2}
\mathcal{N}(\Delta) = \frac{1}{a^{2 - \kappa}}\int_{\Delta}N(x)dx(1 + o(1)), \quad a \to 0,
\end{equation}
where $\mathcal{N}(\Delta)$ is the number of small particles in an arbitrary open subset $\Delta \in D$, $\kappa \in [0,1)$ is a number that can be chosen by an experimenter as desired, $a = \displaystyle\frac{1}{2}\max_{1 \leq m \leq M}\max_{x,y \in D_m}{|x-y|}$ and $N(x) \geq 0$ is a continuous in $D$ function that can be chosen by an experimenter as desired. As $a \to 0$, it is proved (see \cite{[R624]}, \cite{[R635]}) that the solution $u(x,t,a)$ to the problem
\begin{equation}\label{eqi4}
u_t = \Delta u \quad \text{in } D \setminus \displaystyle\bigcup_{m = 1}^M D_m, \, u_N = \zeta_m u \quad \text{on } S_m, \, 1 \leq m \leq M,
\end{equation}
\begin{equation}\label{eqi5}
u|_{t = 0} = f(x), \, \zeta_m := \frac{h(x_m)}{a^\kappa}, \, \text{Re}\zeta_m \geq 0,
\end{equation}
where $x_m \in \mathbb{R}^3$ is an arbitrary point in $D_m$, has a limit, $u(x,t) = \lim_{a \to 0} u(x,t,a)$ which solves problem \eqref{eqi1} with
\begin{equation}\label{eqi6}
q(x) = c_S N(x) h(x), \quad c_S := \frac{|S_m|}{a^2}=const.
\end{equation}

\noindent Therefore, given a potential $q(x)$ (which makes heat propagate along a line), one can choose an arbitrary continuous function $N(x) \geq 0$, can construct a continuous function $h(x) = \frac{q(x)}{c_S N(x)}$, and can distribute the small bodies according to the rule \eqref{eqi2} to obtain the medium in which heat propagates along a line.\\

\noindent Suppose that
\begin{equation}\label{eqi7}
(-\Delta + q(x))\phi(x) = \lambda_n\phi_n, \quad \phi_n|_S = 0, \quad ||\phi_n||_{L^2(D)} = ||\phi_n|| = 1,
\end{equation}
where $\{\phi_n\}$ is an orthonormal basis of $L^2(D) := H$, and $(f, \phi_n):=\int_Df(x)\overline{\phi_n(x)}dx$. Then the unique solution to \eqref{eqi1} is
\begin{equation}\label{eqi8}
u(x, t) = \sum_{n = 1}^\infty e^{-\lambda_n t}(f, \phi_n)\phi_n(x).
\end{equation}
If $q(x)$ is such that $\lambda_1 = 0$, $\lambda_2 \gg 1$, and $\lambda_2 \leq \lambda_3 \leq \dots$,  then, as $t \to \infty$, the series \eqref{eqi8} is well approximated by its first term (see \cite{[R660]} ):
\begin{equation}\label{eqi9}
u(x, t) = (f, \phi_1)\phi_1 + O(e^{-10t}), \quad t \to \infty.
\end{equation}

\noindent Thus, our problem is solved if $q(x)$ has the following property:
\begin{equation}\label{eqi10}
|\phi_1(x)| \text{ decays as $\rho$ grows}, \quad \rho = (x_2^2 + x_3^2)^{1/2}.
\end{equation}
Since the eigenfunction is normalized, $||\phi_1||=1$, this function will not tend to zero in a neighborhood of the line $\rho=0$, so information can be transformed by the heat signals along the line
$\rho=0$, that is, along $s-$axis. Here we use the cylindrical coordinates: $x=(x_1,x_2,x_3)=(s, \rho, \theta)$,
$s=x_1$, $\rho= (x_2^2+x_3^2)^{1/2}$. In Section 2 the potential $q(x)$ will not depend on $\theta$.\\

\noindent The technical part of the solution is the construction of $q(x) = c_S N(x) h(x)$ such that
\begin{equation}\label{eqi11}
\lambda_1 = 0, \quad \lambda_2 \gg 1; \quad |\phi_1(x)| \text{ decays as $\rho$ grows}.
\end{equation}
Since the function $N(x) \geq 0$ and $h(x), \text{Re}h \geq 0$ are at our disposal, any desirable $q, \text{Re}\,q \geq 0$, can be obtained by embedding many small impedance particles in a given domain $D$.

\noindent In section 2, the method for finding such a potential $q(x)$ is presented. In section 3, the numerical method is presented for finding this $q$, and in section 4, the numerical results are presented. In section 5, another numerical method is presented based on the method in section 4.

\section{Construction of $q(x)$ theory}\label{Section2}
\noindent Let $q(x) = p(\rho)+Q(s)$, where $s := x_1$, $\rho:=(x_2^2+x_3^2)^{1/2}$. Then the solution to \eqref{eqi1} is $u = v(\rho)w(s)$, where
\begin{multline}\label{eq1}
-v_m'' - \rho^{-1} v_m' + p(\rho)v_m = \mu_mv_m, \, 0 \leq \rho \leq R, \, |v_m(0)| < \infty, \, v_m(R) = 0;
\end{multline}
\begin{multline}\label{eq2}
-w_j'' + Q(s)w_j = \nu_j w_j, \quad 0 \leq s \leq L, \, w_j(0) = w_j(L) = 0.
\end{multline}

\noindent Our task is to find $Q(s)$ such that $\nu_1 = 0$, $\nu_2 \gg  1$ and $p(\rho)$ such that $\mu_1 = 0, \mu_2 \gg 1$, and $|v_m(\rho)|$ decays as $\rho$ grows.\\

\noindent We use the solution to inverse spectral problem for finding a potential which have the desired properties. The spectral function $\varrho(\lambda)$ of the Dirichlet differential operator $\ell w = -w'' + Q(s)w$ (see formula \eqref{eq2}) is defined by the formula
\[ \varrho(\lambda) = \sum_{\nu_j < \lambda} \frac{1}{\alpha_j}, \]
where $\alpha_j$ are normalizing constants. If $Q = 0$, then the eigenvalues of the corresponding operator $\ell $ are $\nu_{j0} = \displaystyle\left(\frac{\pi j}{L} \right)^2$, $j = 1,2, \dots$, the corresponding normalized eigenfunctions are $\displaystyle\sqrt{\frac{2}{L}}\sin\frac{j \pi x}{L}$, and the normalizing constants $\alpha_{j0} = \displaystyle\sqrt{\int_0^L \sin^2\frac{2\pi x}{L}dx} = \sqrt{\frac{L}{2}}$. If $Q$ is unknown, then the corresponding eigenfunctions are unknown and the role of the normalizing constants can play arbitrary positive numbers which have the right asymptotic. If $L=\pi$ then $\alpha_j=\sqrt{\frac{\pi}{2}} j\big(1+O(\frac 1 j)\big)$, $\sqrt{\nu_j}=j+\frac{c_1}{j} +O(j^{-2})$, and $w_j=\frac{\sin (jx)}{j} +O(j^{-2})$  as $j\to \infty$. \\

Let us recall the procedure, due to Gel'fand and Levitan (see \cite{[L]}, \cite{[R470]}) for finding $Q$ from the
known spectral function. One defines the kernel
$$
L(x,y)=\int_{-\infty}^{\infty} \frac{\sin(\sqrt{\lambda}x)}{\sqrt{\lambda}} \frac{\sin(\sqrt{\lambda}y)}{\sqrt{\lambda}}d(\varrho (\lambda)-\varrho_0(\lambda)),
$$
where $\varrho(\lambda)$ is the spectral function of the operator $\ell $ with the potential $Q=Q(s)$, and  $\varrho_0(\lambda)$ is the spectral function of the operator $\ell $ with the potential $Q=0$ and the same boundary conditions.\\
Consider the Gel'fand-Levitan (GL) integral equation for the kernel $K(x,y)$:
$$ K(x,y)+ \int_0^x K(x,s)L(s,y)ds=-L(x,y), \qquad 0\le y\le x.$$
The solution to this equation allows one to calculate the potential:
$$ Q(x)=2 \frac{d K(x, x)}{dx}.$$

\noindent From now on  we set $L = \pi$.  Then  $\nu_{l0} = j^2$. Let $\nu_1 = 0, \nu_2 = 11, \nu_3 = 14,  \nu_j = \nu_{j0}$ for $j \ge 4$.
Then the kernel $L(x,y)$ in the GL equation is defined as follows
\begin{multline}\label{eq3}
L(x,y) =\frac{3xy}{\pi^3}+ \frac{2}{\pi}\Big(\frac{\sin (\sqrt{\nu_2} x)}{\sqrt{\nu_2}} \frac{\sin (\sqrt{\nu_2} y)}{\sqrt{\nu_2}} +\frac{\sin (\sqrt{\nu_3} x)}{\sqrt{\nu_3}} \frac{\sin (\sqrt{\nu_3} y)}{\sqrt{\nu_3}} \Big) - \\
  -\frac{2}{\pi}\Big(\sin x \sin  y +\sin (2 x) \sin (2 y)  +\sin (3 x) \sin (3 y)\Big),
\end{multline}
where we set the normalizing constants $\alpha_j = \frac{\pi}{2}, j \geq 2, \alpha_1 = \frac{\pi^3}{3}$. The term $xy$
is the value of the function $\frac {\sin \nu x}{\nu} \frac{\sin \nu y}{\nu}$ at $\nu=0$, and $\frac{\pi^3}{3}=||x||^2=\int_0^{\pi}x^2dx$.\\

\noindent Solve the GL equation:
\begin{equation}\label{eq4}
K(s, \tau) + \int_0^sK(s,s')L(s', \tau)ds' = -L(s, \tau), \quad 0 \leq \tau \leq s,
\end{equation}
which is uniquely solvable (see [1], [2]). Equation \eqref{eq4} has finite-rank kernel and therefore can be solved analytically, being equivalent to a linear algebraic system. If $K(s, \tau)$ is found, then
\begin{equation}\label{eq5}
Q(s) = 2\frac{dK(s,s)}{ds},
\end{equation}
and this $Q(s)$ has the required properties: $\nu_1 = 0, \nu_2 \gg 1, \nu_j \leq \nu_{j + 1}$.\\

\noindent Consider now the operator \eqref{eq1} for $v(\rho)$. We want to calculate $p(\rho)$ such that $\mu_1 = 0, \mu_2 \gg 1, \mu_m \leq \mu_{m + 1}, |v_m(\rho)|$ decays as $\rho$ grows.\\

\noindent We reduce this problem to the previous one that was solved. To do this, set $v = \frac{\psi}{\sqrt{\rho}}$. Then equation $-v'' - \frac{1}{\rho}v' + p(\rho)v = \mu v$, is transformed to the equation
\begin{equation}\label{eq6}
-\psi'' - \frac{1}{4\rho^2}\psi + p(\rho)\psi = \mu\psi.
\end{equation}

\noindent Let $p(\rho) = \frac{1}{4\rho^2} + Q(\rho)$, where $Q(\rho)$ is constructed above. Then equation \eqref{eq6} becomes
\begin{equation}\label{eq7}
-\psi'' + Q(\rho)\psi = \mu\psi, \quad \psi(R) = 0, \quad \psi(0)=0.
\end{equation}
It has the desired eigenvalues $\mu_1 = 0, \mu_2 \gg 1, \mu_m \leq \mu_{m + 1}$.\\

\noindent The eigenfunction $\phi_1(x) = v_1(\rho)w_1(s)$, where $v_1(\rho) = \frac{\psi_1(\rho)}{\sqrt{\rho}}$, decays as $\rho$ grows, and the eigenvalues $\lambda_n = \mu_m + \nu_l$. Since $\mu_1 = \nu_1 = 0$ one has $\lambda_1 = 0$. Since $\nu_2 = 11, \mu_2 = 11, \lambda_2 = 11 \gg 1$. Thus, the desired potential is constructed: $q(x) = Q(s)+(\frac{1}{4\rho^2} + Q(\rho))$, where $Q(s)$ is given by formula \eqref{eq5}.\\

\noindent This concludes the description of our procedure for the construction of $q$.

\section{Numerical procedure}\label{Procedure}
In section \ref{NumCons}, a numerical method to construct $q(x)$ is presented. In section \ref{NumCheck}, a procedure is presented to check whether the constructed potential $q(x)$ is valid, by finding the eigenvalues of $q(x)$.
\subsection{Numerical construction of q(x)}\label{NumCons}
From the construction of $q(x)$ in Section \ref{Section2}, if one can construct $Q(s)$ then one gets $q(x) = Q(s) + \left(\frac{1}{4\rho^2} + Q(\rho) \right)$.\\

To construct $Q(s)$, one can use equation \eqref{eq5} and rewrite it as
\begin{equation}\label{eq8}
Q(s) = 2\frac{dK(s,s)}{ds} = 2\left(\frac{\partial K(s,\tau)}{\partial s}\bigg|_{\tau = s} + \frac{\partial K(s,\tau)}{\partial \tau}\bigg|_{\tau = s} \right),
\end{equation}
One can get $K_s$ and $K_\tau$ numerically by the following procedure, using the Gel'fand-Levitan equation \eqref{eq4}.\\

The function $L(x,y)$ in \eqref{eq3} can be written as
\begin{equation}\label{eq9}
L(x,y) = \sum_{j = 1}^6a_j(x)b_j(y),
\end{equation}
where $a_1(x) = \displaystyle\frac{3x}{\pi^3}$, $a_2(x) = \displaystyle\frac{2}{\pi}\displaystyle\frac{\sin(\sqrt{\nu_2}x)}{\nu_2}$, $a_3(x) = \displaystyle\frac{2}{\pi}\frac{\sin(\sqrt{\nu_3}x)}{\nu_3}$, $a_4(x) = -\displaystyle\frac{2}{\pi}\sin(x)$, $a_5(x) = -\displaystyle\frac{2}{\pi}\sin(2x)$, $a_6(x) = -\displaystyle\frac{2}{\pi}\sin(3x)$ and $b_1(x) = x$, $b_2(x) = \sin(\sqrt{\nu_2}x)$, $b_3(x) = \sin(\sqrt{\nu_3}x)$, $b_4(x) = \sin(x)$, $b_5(x) = \sin(2x)$, $b_6(x) = \sin(3x)$.\\

\noindent Then equation \eqref{eq4} becomes
\begin{equation}\label{eq10}
K(s,\tau) + \sum_{j=1}^6b_j(\tau)\int_0^sK(s,s')a_j(s')ds' = -\sum_{j=1}^6a_j(s)b_j(\tau), \quad 0 \leq \tau \leq s.
\end{equation}
Let $\psi_j(s) := \displaystyle\int_0^sK(s,s')a_j(s')ds'$, then equation \eqref{eq10} becomes
\begin{equation}\label{eq11}
K(s,\tau) + \sum_{j=1}^6b_j(\tau)\psi_j(s) = -\sum_{j=1}^6a_j(s)b_j(\tau), \quad 0 \leq \tau \leq s.
\end{equation}
Multiply \eqref{eq11} with $a_m(\tau), 1 \leq m \leq 6$ and integrate it with respect to $\tau$ to get
\begin{equation}\label{eq12}
\psi_m(s) + \sum_{j=1}^6\left( \int_0^sb_j(\tau)a_m(\tau)d\tau \right) \psi_j(s) = -\sum_{j=1}^6a_j(s)\left(\int_0^sb_j(\tau)a_m(\tau)d\tau\right).
\end{equation}
For a fix $s = s_0$, equation \eqref{eq12} is a $6 \times 6$ linear system which can be solved for $\psi_j(s_0), 1 \leq j \leq 6$. So, we can solve equation \eqref{eq12} to get $\psi_j(s)$ numerically. Differentiate equation \eqref{eq12}, one can get a similar linear system to \eqref{eq12} and can get $\psi'_j(s)$ numerically:
\begin{multline}\label{eq13}
\psi'_m(s) + \sum_{j=1}^6\left( \int_0^sb_j(\tau)a_m(\tau)d\tau \right) \psi'_j(s) =\\= -\sum_{j=1}^6\left( a_j(s)\left(\int_0^sb_j(\tau)a_m(\tau)d\tau\right) + (a_j(s) + \psi_j(s))b_j(s)a_m(s) \right).
\end{multline}
After finding $\psi_j(s)$ and $\psi'_j(s)$, one can find $K_s$ and $K_\tau$ by differentiating equation \eqref{eq10} with respect to $s$ and $\tau$
\begin{equation}\label{eq14}
K_s(s,\tau) = -\sum_{j=1}^6\left(a'_j(s)b_j(\tau) + b_j(\tau)\psi'_j(s) \right),
\end{equation}
\begin{equation}\label{eq15}
K_\tau(s,\tau) = -\sum_{j=1}^6\left(a_j(s)b'_j(\tau) + b'_j(\tau)\psi_j(s) \right).
\end{equation}
From equation \eqref{eq14} - \eqref{eq15}, one finds $K_s(s,s)$ and $K_\tau(s,s)$, and then finds $Q(s)$ numerically using equation \eqref{eq8}.

\subsection{Checking the eigenvalues of q(x)}\label{NumCheck}
To check whether the constructed potential $q(x)$ is the correct potential, one has to check whether the eigenvalues generated by the constructed potential $q(x)$ satisfy the conditions in Section \ref{Section2}. It is sufficient to check the eigenvalues of $Q(s)$ are $\nu_1 = 0, \nu_2 = 11, \nu_3 = 14, \nu_j = j^2, j \geq 4$.\\

One can find numerically the eigenvalues of $Q(s)$ using the variational definition of eigenvalues
\begin{equation}\label{eq16}
\nu_j = \inf_{u \neq 0} \frac{\int_0^\pi (|u'|^2 + Q|u|^2)dx}{\int_0^\pi|u|^2dx}, \quad u = \sum_{n = 1}^\infty c_n \varphi_n(x),
\end{equation}
where $\varphi_n(x) := \sqrt{\frac{2}{\pi}}\sin(nx), n \geq 1$ is the normalized eigenfunctions of
\begin{equation}\label{eq17}
-\varphi'' + Q(s)\varphi = \nu\varphi, \quad \varphi(0) = \varphi(\pi) = 0, \quad 0 \leq s \leq \pi.
\end{equation}
Finding $\nu_j$ by equation \eqref{eq16} is equivalent to finding the infimum of
\begin{equation}\label{eq18}
\sum_{n = 1}^N n^2|c_n|^2 + \sum_{n,m = 1}^Nc_n\bar{c_m}q_{nm}, \text{ under the restriction } \displaystyle\sum_{n=1}^N|c_n|^2 = 1,
\end{equation}
where $q_{nm} := \displaystyle\frac{2}{\pi}\int_0^\pi Q(x)\sin(nx)\sin(mx)dx$. Minimizing \eqref{eq18} is equivalent to minimizing
\begin{equation}\label{eq19}
f(c_1, \dots, c_N) = \sum_{n = 1}^N n^2|c_n|^2 + \sum_{n,m = 1}^Nc_n\bar{c_m}q_{nm} - \nu\left(\sum_{n = 1}^N c_n^2 - 1\right).
\end{equation}
The function $f(c_1, \dots, c_N)$ attains the minimum when $\displaystyle\frac{\partial f}{\partial c_n} = 0$. This leads to a linear system with respect to $c_1, \dots, c_N$:
\begin{equation}\label{eq20}
n^2c_n + \sum_{m=1}^Nq_{nm}c_m - \nu c_n = 0.
\end{equation}
This linear system can be written as
\begin{equation}\label{eq21}
PC = \nu C,
\end{equation}
where $P$ is a matrix such that $P_{nm} = n^2\delta_{nm} + q_{nm}$, and $C$ is a column vector $C = (c_1, \dots, c_N)$. Then the eigenvalues of $Q(s)$ are the eigenvalues of the matrix $P$.

\subsection{Calculating $q_{nm}$}\label{Calqnm}
\noindent In section \ref{NumCheck}, one needs to calculate $q_{nm} = \displaystyle\frac{2}{\pi}\int_0^\pi Q(x)\sin(nx)\sin(mx)dx$ to construct the matrix $P$. One can calculate the matrix $P$ by using the formula $2\sin A \sin B = \cos(A - B) - \cos(A + B)$. So, one first calculates $$\tilde{q}(k) := \displaystyle\frac{1}{\pi}\int_0^\pi Q(x) \cos(kx)dx, 0 \leq k \leq 2N,$$ and calculates $q_{nm}$ by the formula $$q_{nm} = \tilde{q}(|n-m|) - \tilde{q}(n+m)$$.

\section{Numerical results}\label{NumResults}
\noindent Based on the numerical procedure in section \ref{Procedure}, a computer algorithm/program is developed with the following main steps:
\begin{itemize}
\item[1.] Partition the interval $[0, \pi]$ into $\mathcal{M}$ equal-distanced intervals with the endpoints $x_i$, $1 \leq i \leq \mathcal{M} + 1$.
\item[2.] For each $x_i, 1 \leq i \leq \mathcal{M} + 1$, one solves the linear systems \eqref{eq12} and \eqref{eq13} for $\psi_j(x_i)$ and $\psi'_j(x_i)$, $1 \leq j \leq 6$.
\item[3.] Find $K_s(x_i,x_i)$ and $K_\tau(x_i, x_i)$ by using equations \eqref{eq14} and \eqref{eq15}.
\item[4.] Find $Q(x_i), 1 \leq i \leq \mathcal{M} + 1$.
\item[5.] Construct the matrix $P$ in equation \eqref{eq21} by calculating $q_{nm}$ using the procedure in section 3.3, where $\tilde{q}(k)$ is calculated using the Riemann sum
$$\tilde{q}(k) = \frac{1}{\pi}\sum_{i = 1}^\mathcal{M} \frac{Q(x_i)\cos(kx_i) + Q(x_{i+1})\cos(kx_{i+1})}{2}(x_{i + 1} - x_i).$$
\item[6.] Find the eigenvalues of $P$ using the Jacobi eigenvalue algorithm see, for example \cite{[Jacobi]}.
\end{itemize}

\noindent The above algorithm is run five times with $\mathcal{M} = 100, 150, 200, 250$, and $300$. The constructed potentials $Q(s)$ are as in the following graph.

\begin{figure}[H]
  \centering
      \includegraphics[width=300px]{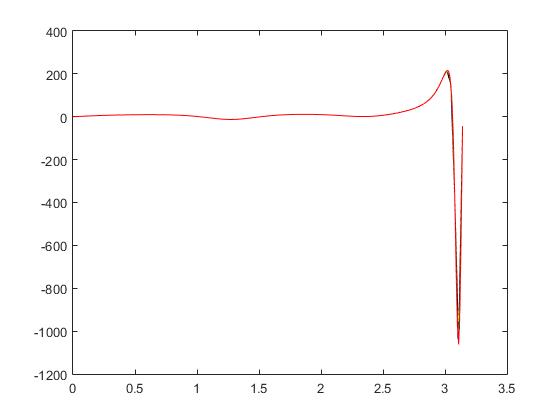}
  \caption{Numerically constructed potentials $Q(s)$.}
\end{figure}
\noindent The relative error of the eigenvalues is calculated by the following formula:
$$\delta_\mathcal{M} = \max_{1 \leq j \leq \mathcal{M}}\frac{|\nu_j^{(\mathcal{M})} - \nu_j|}{\nu_j}.$$
The following table gives the relative errors of the eigenvalues of the constructed potentials $Q(s)$ for $\mathcal{M} = 100, 150, 200, 250$, and $300$.
\\
\begin{center}
	\begin{tabular}{| l | c |}
    \hline
    $\mathcal{M}$ & Relative errors $\delta_\mathcal{M}$ \\ \hline
    100 & 57.12\% \\ \hline
    150 & 5.01\% \\ \hline
    200 & 1.65\% \\ \hline
    250 & 0.67\% \\ \hline
    300 & 0.32\% \\ \hline
    \end{tabular}
\end{center}
From the above table, one can construct the potential $Q(s)$ with $\mathcal{M} = 250$ equal-distance small intervals and gets the relative error of the eigenvalues less than $1\%$. The above result also shows that the constructed $Q(s)$ is valid.

\section{Another method to calculate eigenvalues}\label{Sec5}
In the numerical results in section \ref{NumResults}, one needs $\mathcal{M} = 250$ equal-distance small intervals to get the relative error less than $1\%$. From the graph of the constructed $Q(s)$, since $Q(s)$ is pretty steep close to $\pi$, one can improve the method in section \ref{NumResults} by distributing more equal-distanced intervals $\mathcal{M}_2$ in the interval $[\frac{9\pi}{10},\pi]$ and less equal-distanced intervals $\mathcal{M}_1$ in the interval $[0, \frac{9\pi}{10}]$.\\

The following result is obtained
\begin{figure}[H]
  \centering
      \includegraphics[width=300px]{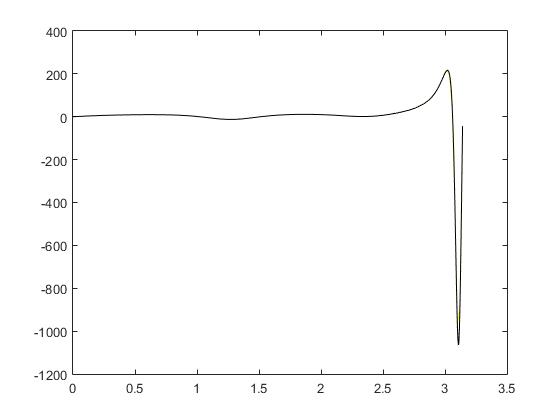}
  \caption{Numerically constructed potential $Q(s)$ for non-uniform small intervals.}
\end{figure}

\begin{center}
	\begin{tabular}{| l | l | c |}
    \hline
    $\mathcal{M}_1$ & $\mathcal{M}_2$ & Relative errors \\ \hline
    50 & 50 & 4.68\% \\ \hline
    50 & 75 & 0.94\% \\ \hline
    50 & 100 & 0.23\% \\ \hline
    \end{tabular}
\end{center}

\begin{remark}
About calculating $\tilde{q}(k)$
\end{remark}
\noindent In the constructed potential $Q(s)$ in Figure 2, let $Q_{max} := \displaystyle\max_{1 \leq i \leq \mathcal{M} + 1} Q(x_i)$, $x_{max} := \{x_i: Q(x_i) = Q_{max} \}$, $Q_{min} := \displaystyle\min_{1 \leq i \leq \mathcal{M} + 1} Q(x_i)$, and $x_{min} := \{x_i: Q(x_i) = Q_{min} \}$. The parts of the potential from $x_{max}$ to $x_{min}$ and from $x_{min}$ to $\pi$ look like straight lines and one may try to calculate $\tilde{q}(k)$ by
\begin{multline}\label{eq51}
\tilde{q}(k) = \displaystyle\frac{1}{\pi} \left( \int_0^{x_{max}} Q(x) \cos(kx)dx \, +\right. \\ \left. + \int_{x_{max}}^{x_{min}} Q(x) \cos(kx)dx + \int_{x_{min}}^{\pi} Q(x) \cos(kx)dx \right)\\
:= I_1 + I_2 + I_3.
\end{multline}
If $Q(x)$ is a straight line from $x_{max}$ to $x_{min}$ and from $x_{min}$ to $\pi$, one can calculate $I_2$ and $I_3$ analytically. However, this does not provide the desired numerical accuracy as the following numerical experiment shows.\\

\noindent In this experiment, $P_1$ is the numerical matrix in equation \eqref{eq21} obtained in the experiment described in Section 5 for $\mathcal{M}_1 = 50$ and $\mathcal{M}_2 = 75$, and $P_2$ is the numerical matrix obtained by considering $Q(x)$ as straight lines from $x_{max}$ to $x_{min}$ and from $x_{min}$ to $\pi$. The relative error matrix $E$ is calculated by $E_{nm} = \displaystyle\frac{|P_{1,nm} - P_{2,nm}|}{|P_{1,nm}|}$. Then $\min(E) = 0.37$ but $\max(E) = 510.17$. So, although the parts of the potential look like straight lines from $x_{max}$ to $x_{min}$ and from $x_{min}$ to $\pi$, one cannot consider them as straight lines in numerical calculations.

\newpage
\noindent\textbf{Conflict of Interest:} The authors declare that they have no conflict of interest.

\end{document}